# REDUCED ORDER DEAD-BEAT OBSERVERS FOR A BIOREACTOR

**Iasson Karafyllis[*] and Zhong-Ping Jiang[**]**

[*]Dept. of Environmental Eng., Technical University of Crete,
73100, Chania, Greece, email: ikarafyl@enveng.tuc.gr

[**]Dept. of Electrical and Computer Eng., Polytechnic Institute of New York University,
Six Metrotech Center, Brooklyn, NY 11201, U.S.A., email: zjiang@control.poly.edu

**Abstract**
This paper studies the strong observability property and the reduced-order dead-beat observer design problem for a continuous bioreactor. New relationships between coexistence and strong observability, and checkable sufficient conditions for strong observability, are established for a chemostat with two competing microbial species. Furthermore, the dynamic output feedback stabilization problem is solved for the case of one species.

**Keywords:** observer design, chemostat, hybrid systems.

## 1. Introduction

The design of observers is a very important problem in mathematical control theory. In this work we focus on the observer design problem for the chemostat with $n$ competing species and one limiting substrate (see [30]):

$$\dot{x}_i(t) = \left(\mu_i(s(t)) - D(t) - b_i\right)x_i(t) \quad , \quad i = 1,...,n$$
$$\dot{s}(t) = D(t)\left(s_{in}(t) - s(t)\right) - \sum_{i=1}^{n} g_i(s(t))x_i(t) \quad (1.1)$$
$$x(t) = (x_1(t),...,x_n(t))' \in \Omega \, , \, s(t) \in A$$

with measurement $y = s$ and inputs $u = (D, s_{in}) \in U = [0,+\infty)^2 \subseteq \Re^2$. As usual $x_i(t)$ denotes the concentration of the $i$th microbial species and $\Omega = \text{int}(\Re_+^n)$, $s(t) \in A \subseteq (0,+\infty)$ denotes the concentration of the limiting nutrient, $s_{in}(t)$ denotes the inlet concentration of the limiting nutrient and $D(t)$ denotes the dilution rate. The specific growth rate $\mu_i(s)$ of the $i$th microbial species is a continuously differentiable, bounded function $\mu_i : [0,+\infty) \to [0,+\infty)$ with $\mu_i(0) = 0$ and $\mu_i(s) > 0$ for all $s > 0$ ($i = 1,...,n$). The constants $b_i \geq 0$ ($i = 1,...,n$) are the mortality rates of the microbial species, while the continuously differentiable, bounded functions $g_i : [0,+\infty) \to [0,+\infty)$ with $g_i(0) = 0$ and $g_i(s) > 0$ for all $s > 0$ ($i = 1,...,n$) are the products of the specific growth rates of the species with the corresponding possibly varying yield constants (see [32,33] for the chemostat with variable yields). For the open set $A \subseteq (0,+\infty)$, we will distinguish the following cases:

- $A = (0,+\infty)$ (the general case),

- for the case $s_{in}(t) \equiv s_{in}$, $A = (0, s_{in})$.

The literature concerning chemostat models of the form (1.1) is vast, since the chemostat appears to be one of the cornerstones of Mathematical Population Biology. The dynamics of (1.1) were studied in [7,25,30,32,33] (see also references therein), where the theory of monotone dynamical systems (see [2,29,30]) plays an important role. Feedback stabilization problems for the chemostat have been studied in [3,6,12,13,14,16,17,18,19,23,24,26,27,28].

Observer problems for the chemostat have been studied in [3,4,5,8,9,10,11,22]. It should be noted that the successful observer designs for the chemostat can lead to the solution of dynamic output feedback control problems (see [3]).

In this work, we apply recent results for the observability of systems which are linear to unmeasured state components in [20] in order to design hybrid dead-beat reduced-order observers for system (1.1). More specifically, we show that:

- system (1.1) with $n=1$ is strongly observable in time $r>0$ for arbitrary $r>0$,

- the dynamic output feedback stabilization problem for (1.1) with $n=1$ can be solved with the combination of a static state feedback stabilizer and the proposed dead-beat hybrid reduced order observer (Proposition 3.1),

- coexistence implies absence of strong observability for $n=2$ (Proposition 3.3),

- it is impossible to design a smooth observer for system (1.1) with $n=2$ which guarantees convergence of the estimates for all inputs $(D, s_{in}) \in L^{\infty}_{loc}(\Re_+; U)$, under the assumption of coexistence,

- system (1.1) with $n=2$, $D \equiv 0$ (batch culture) and Michaelis-Menten kinetics for the specific growth rates is strongly observable in time $r>0$ for arbitrary $r>0$ if and only if the specific growth rates and the mortality rates of the two microbial species are not identical (Theorem 3.5); this result is important because batch cultures of microbial species are used only for finite time and the proposed dead-beat hybrid reduced order observer can provide exact estimates in very short times,

- a set of sufficient conditions (which does not allow coexistence) can guarantee strong observability of system (1.1) with $n=2$ in time $r>0$ (Theorem 3.8 and Theorem 3.9).

**Notations** Throughout this paper we adopt the following notations:

* Let $I \subseteq \Re_+ := [0, +\infty)$ be an interval. By $L^{\infty}(I; U)$ ($L^{\infty}_{loc}(I; U)$) we denote the space of measurable and (locally) essentially bounded functions $u(\cdot)$ defined on $I$ and taking values in $U \subseteq \Re^m$.
* For definitions of the function classes $K_{\infty}$, $KL$ see [21].
* By $C^0(A; \Omega)$, we denote the class of continuous functions on $A$, which take values in $\Omega$.
* For a vector $x \in \Re^n$ we denote by $x'$ its transpose. The determinant of a square matrix $A \in \Re^{n \times n}$ is denoted by $\det(A)$. $A' \in \Re^{n \times m}$ denotes the transpose of the matrix $A \in \Re^{m \times n}$.
* By $A = diag(l_1, l_2, ..., l_n)$ we mean that the matrix $A = \{a_{ij}; i=1,...,n, j=1,...,n\}$ is diagonal with $a_{ii} = l_i$, for $i=1,...,n$.
* By $int(U) \subseteq U$ we denote the interior of a set $U \subseteq \Re^m$. By $\Re^n_+$ we denote the set of all $x = (x_1, ..., x_n)' \in \Re^n$ with $x_i \geq 0$ ($i=1,...,n$).

Under the assumption that the continuously differentiable functions $\mu_i : [0, +\infty) \to [0, +\infty)$ with $\mu_i(0) = 0$ and $\mu_i(s) > 0$ for all $s > 0$, $g_i : [0, +\infty) \to [0, +\infty)$ with $g_i(0) = 0$ and $g_i(s) > 0$ for all $s > 0$ ($i=1,...,n$) are bounded, it follows that for every pair of inputs $(D, s_{in}) \in L^{\infty}_{loc}(\Re_+; U)$ and for every $(x_0, s_0) \in \Omega \times A$ there exists a unique solution $(x(t), s(t)) \in \Omega \times A$ of (1.1) defined for all $t \geq 0$ with initial condition $(x(0), s(0)) = (x_0, s_0)$ corresponding to inputs $(D, s_{in}) \in L^{\infty}_{loc}(\Re_+; U)$.

## 2. Review of Recent Results

Consider an autonomous system described by ordinary differential equations of the form:

$$\dot{x}(t) = f(x(t), u(t))$$
$$x(t) \in \Omega \subseteq \Re^n, u(t) \in U \subseteq \Re^m \tag{2.1}$$

where $\Omega \subseteq \Re^n$ is an open set, $U \subseteq \Re^m$ is a non-empty closed set and the mapping $f: \Omega \times U \to \Re^n$ is locally Lipschitz. The output of system (2.1) is given by

$$y(t) = h(x(t)) \tag{2.2}$$

where the mapping $h: \Omega \to \Re^k$ is continuous. We assume that for every $x_0 \in \Omega$ and $u \in L^\infty_{loc}(\Re_+; U)$, the solution $x(t, x_0; u)$ of (2.1) with initial condition $x(0) = x(0, x_0; u) = x_0$ and corresponding to input $u \in L^\infty_{loc}(\Re_+; U)$ exists for all $t \geq 0$, i.e., we assume forward completeness. For system (2.1) we adopt the following notion of observability.

**Definition 2.1 (Definition 2.1 in [20]):** *Consider system (2.1) with output (2.2). We say that the input $u \in L^\infty([0, r]; U)$ strongly distinguishes the state $x_0 \in \Omega$ in time $r > 0$, if the following condition holds*

$$\max_{t \in [0, r]} |h(x(t, x_0; u)) - h(x(t, \xi; u))| > 0, \text{ for all } \xi \in D \text{ with } x_0 \neq \xi \tag{2.3}$$

We next define the notion of strongly observable systems in time $r > 0$. The reader should note that strong observability is a more demanding notion than simple observability (see [31]).

**Definition 2.2 (Definition 2.7 in [20]):** *Consider system (2.1). We say that (2.1) is strongly observable in time $r > 0$ if every input $u \in L^\infty([0, r]; U)$ strongly distinguishes every state $x_0 \in D$ in time $r > 0$.*

Now we show how the main results of [20], that is, Proposition 2.3 and Corollary 2.11 of [20], can be applied to the continuous bioreactor model (1.1). More precisely, assuming for the time being that system (1.1) is strongly observable in time $r > 0$, we are in a position to define the operator:

$$P: C^0([0, r]; A) \times L^\infty([0, r]; U) \to \Omega$$

For each $s \in C^0([0, r]; A)$, $(D, s_{in}) \in L^\infty([0, r]; U)$, $P(s, D, s_{in})$ is defined by

$$P(s, D, s_{in}) = \Phi(r) Q^{-1} \int_0^r p(\tau) q(\tau) d\tau \tag{2.4}$$

where $\Phi(t) := diag\left(\exp\left(\int_0^t (\mu_1(s(w)) - D(w) - b_1) dw\right), \ldots, \exp\left(\int_0^t (\mu_n(s(w)) - D(w) - b_n) dw\right)\right)$, $Q = \int_0^r q(\tau) q'(\tau) d\tau$, $q(\tau) = \int_0^\tau \Phi'(s) C(s) ds$, $C(\tau) = -(g_1(s(\tau)), \ldots, g_n(s(\tau)))' \in \Re^n$, $p(\tau) = s(\tau) - s(0) - \int_0^\tau D(w)(s_{in}(w) - s(w)) dw$ for all $\tau \in [0, r]$. Proposition 2.3 in [20] guarantees that, under the assumption of strong observability in time $r > 0$ for system (1.1), then for every $(x_0, s_0) \in \Omega \times A$ and $(D, s_{in}) \in L^\infty_{loc}(\Re_+; U)$ the following equality holds:

$$x(t) = P(\delta_{t-r} s, \delta_{t-r} D, \delta_{t-r} s_{in}), \text{ for all } t \geq r \tag{2.5}$$

where $(\delta_{t-r} s)(w) = s(t - r + w)$, $(\delta_{t-r} D)(w) = D(t - r + w)$, $(\delta_{t-r} s_{in})(w) = s_{in}(t - r + w)$ for $w \in [0, r]$.

Therefore, if system (1.1) is strongly observable in time $r > 0$, then we are in a position to provide a reduced order dead-beat observer for system (1.1). Given $t_0 \geq 0$, $z_0 \in \Omega$, we calculate $z(t)$ by the following algorithm:

Step $i$: Calculation of $z(t)$ for $t \in [t_0 + ir, t_0 + (i+1)r]$

1) Calculate $z(t)$ for $t \in [t_0 + ir, t_0 + (i+1)r)$, the solution of $\dot{z}(t) = M(s(t), D(t)) z(t)$, where $M(s, D) = diag(\mu_1(s) - D - b_1, \ldots, \mu_n(s) - D - b_n)$.

2) Set $z(t_0 + (i+1)r) = P(\delta_{t_0 + ir} s, \delta_{t_0 + ir} D, \delta_{t_0 + ir} s_{in})$, where $P: C^0([0, r]; A) \times L^\infty([0, r]; U) \to \Omega$ is the operator defined by (2.4).

For $i = 0$ we take $z(t_0) = z_0$ (initial condition). Schematically, we write:

$$\begin{aligned}
\dot{z}(t) &= M(s(t), D(t))z(t), \, t \in [\tau_i, \tau_{i+1}) \\
z(\tau_{i+1}) &= P(\delta_{\tau_i} s, \delta_{\tau_i} D, \delta_{\tau_i} s_{in}) \\
\tau_{i+1} &= \tau_i + r
\end{aligned} \tag{2.6}$$

Thus, we obtain from Corollary 2.11 in [20] the following result.

**Proposition 2.3:** *Consider system (1.1) and assume that it is strongly observable in time $r > 0$. Consider the unique solution $(x(t), s(t), z(t)) \in \Omega \times A \times \Omega$ of (1.1), (2.6) with arbitrary initial condition $(x_0, s_0, z_0) \in \Omega \times A \times \Omega$ corresponding to arbitrary input $(D, s_{in}) \in L^\infty_{loc}(\Re_+; U)$. Then the solution $(x(t), s(t), z(t)) \in \Omega \times A \times \Omega$ of (1.1), (2.6) satisfies:*

$$z(t) = x(t), \text{ for all } t \geq r \tag{2.7}$$

The explicit formulae for the observer (2.6) for $n = 2$ are given next.

$$\begin{aligned}
\dot{z}_1(t) &= \left(\mu_1(s(t)) - D(t) - b_1\right) z_1(t) \\
\dot{z}_2(t) &= \left(\mu_2(s(t)) - D(t) - b_2\right) z_2(t)
\end{aligned}, \text{ for } t \in [\tau_i, \tau_{i+1})$$

and

$$z_1(\tau_{i+1}) = \frac{N_1 \exp\left(\int_{\tau_i}^{\tau_{i+1}} (\mu_1(s(t)) - D(t) - b_1) dt\right)}{\int_{\tau_i}^{\tau_{i+1}} \phi_1^2(t) dt \int_{\tau_i}^{\tau_{i+1}} \phi_2^2(t) dt - \left(\int_{\tau_i}^{\tau_{i+1}} \phi_1(t) \phi_2(t) dt\right)^2}; \; z_2(\tau_{i+1}) = \frac{N_2 \exp\left(\int_{\tau_i}^{\tau_{i+1}} (\mu_{21}(s(t)) - D(t) - b_2) dt\right)}{\int_{\tau_i}^{\tau_{i+1}} \phi_1^2(t) dt \int_{\tau_i}^{\tau_{i+1}} \phi_2^2(t) dt - \left(\int_{\tau_i}^{\tau_{i+1}} \phi_1(t) \phi_2(t) dt\right)^2}$$

where

$$\begin{aligned}
N_1 &:= \int_{\tau_i}^{\tau_{i+1}} \phi_2^2(t) dt \int_{\tau_i}^{\tau_{i+1}} \left( s(\tau_i) - s(t) + \int_{\tau_i}^{t} D(w)(s_{in}(w) - s(w)) dw \right) \phi_1(t) dt \\
&\quad - \int_{\tau_i}^{\tau_{i+1}} \phi_1(t) \phi_2(t) dt \int_{\tau_i}^{\tau_{i+1}} \left( s(\tau_i) - s(t) + \int_{\tau_i}^{t} D(w)(s_{in}(w) - s(w)) dw \right) \phi_2(t) dt \\
N_2 &:= \int_{\tau_i}^{\tau_{i+1}} \phi_1^2(t) dt \int_{\tau_i}^{\tau_{i+1}} \left( s(\tau_i) - s(t) + \int_{\tau_i}^{t} D(w)(s_{in}(w) - s(w)) dw \right) \phi_2(t) dt \\
&\quad - \int_{\tau_i}^{\tau_{i+1}} \phi_1(t) \phi_2(t) dt \int_{\tau_i}^{\tau_{i+1}} \left( s(\tau_i) - s(t) + \int_{\tau_i}^{t} D(w)(s_{in}(w) - s(w)) dw \right) \phi_1(t) dt
\end{aligned}$$

and

$$\begin{aligned}
\phi_1(t) &:= \int_{\tau_i}^{t} g_1(s(\tau)) \exp\left( \int_{\tau_i}^{\tau} (\mu_1(s(w)) - D(w) - b_1) dw \right) d\tau \\
\phi_2(t) &:= \int_{\tau_i}^{t} g_2(s(\tau)) \exp\left( \int_{\tau_i}^{\tau} (\mu_2(s(w)) - D(w) - b_2) dw \right) d\tau
\end{aligned}$$

It should be noted that the observer (2.6) is a hybrid observer which uses delays and guarantees exact knowledge of the concentrations of the microbial species after $r$ time units.

Finally, we end this section by presenting the following lemma, which is to be used in next section.

**Lemma 2.4:** *Suppose that system (1.1) with $n=2$ is not strongly observable in time $r>0$. Then there exist $(D, s_{in}) \in L^\infty([0,r];U)$ and $(x_0, s_0) \in \text{int}(\Re_+^2) \times A$ such that*

$$\dot{s}(t) = D(t)(s_{in}(t) - s(t)) - \left(x_{2,0} + \frac{g_1(s_0)}{g_2(s_0)} x_{1,0}\right) g_2(s(t)) \exp\left(\int_0^t (\mu_2(s(w)) - D(w) - b_2) dw\right), \text{ for almost all } t \in [0,r] \quad (2.8)$$

*and*

$$\kappa(s(t))\dot{s}(t) = \mu_2(s(t)) - \mu_1(s(t)) + b_1 - b_2, \text{ for almost all } t \in [0,r] \quad (2.9)$$

*where $s(t) \in A$ denotes the component of the solution $(x(t), s(t)) \in \text{int}(\Re_+^2) \times A$ of (1.1) with initial condition $(x_0, s_0) \in \text{int}(\Re_+^2) \times A$ corresponding to input $(D, s_{in}) \in L^\infty_{loc}(\Re_+; U)$, $x_0 = (x_{1,0}, x_{2,0}) \in \text{int}(\Re_+^2)$ and*

$$\kappa(s) := \frac{d}{ds} \ln\left(\frac{g_1(s)}{g_2(s)}\right) \quad (2.10)$$

**Proof:** Suppose that system (1.1) is not strongly observable in time $r>0$. By virtue of Remark 2.4 in [20] there exists $(D, s_{in}) \in L^\infty([0,r];U)$, $(x_0, s_0) \in \text{int}(\Re_+^2) \times A$ and $\xi = (\xi_1, \xi_2)' \in \Re^2$, $\xi \neq 0$ such that

$$\xi_1 g_1(s(t)) \exp\left(\int_0^t (\mu_1(s(w)) - b_1) dw\right) + \xi_2 g_2(s(t)) \exp\left(\int_0^t (\mu_2(s(w)) - b_2) dw\right) = 0, \text{ for all } t \in [0,r] \quad (2.11)$$

where $s(t) \in A$ denotes the component of the solution $(x(t), s(t)) \in \text{int}(\Re_+^2) \times A$ of (1.1) with initial condition $(x_0, s_0) \in \text{int}(\Re_+^3)$ corresponding to input $(D, s_{in}) \in L^\infty_{loc}(\Re_+; U)$. Therefore, (2.11) implies that:

$$\det\left(\begin{bmatrix} g_1(s_0) & g_2(s_0) \\ g_1(s(t)) \exp\left(\int_0^t (\mu_1(s(w)) - b_1) dw\right) & g_2(s(t)) \exp\left(\int_0^t (\mu_2(s(w)) - b_2) dw\right) \end{bmatrix}\right) = 0, \text{ for all } t \in [0,r] \quad (2.12)$$

or equivalently,

$$g_1(s(t)) \exp\left(\int_0^t (\mu_1(s(w)) - b_1) dw\right) = \frac{g_1(s_0)}{g_2(s_0)} g_2(s(t)) \exp\left(\int_0^t (\mu_2(s(w)) - b_2) dw\right), \text{ for all } t \in [0,r] \quad (2.13)$$

Equality (2.13) in conjunction with (1.1) and the fact that $x_i(t) = x_{i,0} \exp\left(\int_0^t (\mu_i(s(w)) - D(w) - b_i) dw\right)$, $i=1,2$, implies equality (2.8). Equality (2.9) is obtained by differentiation of (2.13). ◁

## 3. Strong Observability of the Chemostat

The results of the previous section clearly indicate that it is important to study under what conditions the chemostat (1.1) is a strongly observable system.

### 3.A. Case $n=1$

For this case we have the system:

$$\begin{aligned} \dot{x}(t) &= (\mu(s(t)) - D(t) - b) x(t) \\ \dot{s}(t) &= D(t)(s_{in}(t) - s(t)) - g(s(t)) x(t) \end{aligned} \quad (3.1)$$

Using Corollary 2.5 in [20], we conclude that system (3.1) is strongly observable in time $r>0$ for arbitrary $r>0$. The observer (2.6) can be used for system (3.1), which for $n=1$ takes the form:

$$\dot{z}(t) = (\mu(s(t)) - D(t) - b)z(t), \ t \in [\tau_i, \tau_{i+1})$$

$$z(\tau_{i+1}) = \exp\left(\int_{\tau_i}^{\tau_{i+1}} (\mu(s(w)) - D(w) - b)dw\right) \frac{\int_{\tau_i}^{\tau_{i+1}} p(t)q(t)dt}{\int_{\tau_i}^{\tau_{i+1}} q^2(t)dt} \quad (3.2)$$

$$\tau_{i+1} = \tau_i + r$$

where

$$p(t) = s(t) - s(\tau_i) - \int_{\tau_i}^{t} D(w)(s_{in}(w) - s(w))dw$$

$$q(t) = -\int_{\tau_i}^{t} g(s(p)) \exp\left(\int_{\tau_i}^{p} (\mu(s(w)) - D(w) - b)dw\right)dp$$

The hybrid reduced-order observer (3.2) guarantees that for every initial condition $(x_0, s_0, z_0) \in \Omega \times A \times \Omega$ and for every $(D, s_{in}) \in L^\infty_{loc}(\Re_+; U)$ the solution of (3.1), (3.2) satisfies $z(t) = x(t)$ for all $t \geq r$.

It is important to notice that the observer (3.2) can be combined with state feedback laws for the stabilization of the chemostat model (3.1) by means of dynamic output feedback. More specifically, we obtain the following result.

**Proposition 3.1:** *Let $(s^*, x^*) \in (0, s_{in}) \times (0, +\infty)$ be an equilibrium point for (3.1) with $s_{in}(t) \equiv s_{in} > 0$, $D(t) \equiv D^* > 0$, i.e., $\mu(s^*) = D^* + b$ and $D^*(s_{in} - s^*) = g(s^*)x^*$. Suppose that there exists a locally Lipschitz feedback law $k:(0, s_{in}) \times (0, +\infty) \to (0, +\infty)$ with $D^* = k(s^*, x^*)$ such that $(s^*, x^*) \in (0, s_{in}) \times (0, +\infty)$ is globally asymptotically stable for the closed-loop system (3.1) with the feedback law $D = k(s,x)$ and $(s,x) \in (0, s_{in}) \times (0, +\infty)$. Moreover, suppose that system (3.1) with the dynamic feedback law $D = k(s,z)$, $\dot{z} = (\mu(s) - D - b)z$, $z \in (0, +\infty)$ is forward complete.*

*Then for every $r > 0$, the equilibrium point $(s, x, z) = (s^*, x^*, x^*)$ is globally asymptotically stable for the closed-loop system (3.1) with (3.2) and $D = k(s,z)$.*

**Proof:** Let $r > 0$ be arbitrary. First, the following change of coordinates is performed:

$$x = x^* \exp(x_1), \ s = \frac{s_{in} \exp(x_2)}{G + \exp(x_2)} \quad (3.3)$$

where $G := \frac{s_{in} - s^*}{s^*}$. Under the above change of coordinates, system (3.1) with $s_{in}(t) \equiv s_{in} > 0$ takes the form:

$$\dot{x}_1 = \tilde{\mu}(x_2) + D^* - D$$
$$\dot{x}_2 = (G\exp(-x_2) + 1)\left[D - D^*\tilde{g}(x_2)\exp(x_1)\right] \quad (3.4)$$
$$x = (x_1, x_2)' \in \Re^2$$

where $\tilde{\mu}(x_2) := \mu\left(\frac{s_{in}\exp(x_2)}{G+\exp(x_2)}\right) - \mu(s^*)$, $\tilde{g}(x_2) := \frac{G+\exp(x_2)}{G+1} \frac{g\left(\frac{s_{in}\exp(x_2)}{G+\exp(x_2)}\right)}{g(s^*)}$. Since the feedback law $D = k(s,x)$ guarantees that $(s^*, x^*) \in (0, s_{in}) \times (0, +\infty)$ is globally asymptotically stable for the closed-loop system (3.1), it follows that $0 \in \Re^2$ is globally asymptotically stable for the closed-loop system (3.4) with $D = k\left(\dfrac{s_{in}\exp(x_2)}{G+\exp(x_2)}, x^*\exp(x_1)\right)$. Therefore, there exists $\sigma \in KL$ such that for every $x_0 \in \Re^2$, $t \geq 0$ the solution of the closed-loop system (3.4) with $D = k\left(\dfrac{s_{in}\exp(x_2)}{G+\exp(x_2)}, x^*\exp(x_1)\right)$ initiated from $(x_1(0), x_2(0))' = x_0$ satisfies:

$$|x(t)| \leq \sigma(|x_0|, t) \qquad (3.5)$$

System (3.1) with the dynamic feedback law $D = k(s, z)$, $\dot{z} = (\mu(s) - D - b)z$ is transformed to system (3.4) with $D = k\left(\dfrac{s_{in}\exp(x_2)}{G+\exp(x_2)}, x^*\exp(x_3)\right)$ and

$$\dot{x}_3 = \tilde{\mu}(x_2) + D^* - D \quad ; \quad x_3 \in \Re \qquad (3.6)$$

where

$$z = x^*\exp(x_3) \qquad (3.7)$$

Moreover, system (3.4), (3.6) with $D = k\left(\dfrac{s_{in}\exp(x_2)}{G+\exp(x_2)}, x^*\exp(x_3)\right)$ is forward complete. Using the main result in [1], it follows that system (3.4), (3.6) with $D = k\left(\dfrac{s_{in}\exp(x_2)}{G+\exp(x_2)}, x^*\exp(x_3)\right)$ is Robustly Forward Complete (see [15]) and consequently, by virtue of Lemma 2.3 in [15], there exists $a \in K_\infty$ such that every $x_0 \in \Re^3$, $t \in [0, r]$ the solution of the closed-loop system (3.4), (3.6) with $D = k\left(\dfrac{s_{in}\exp(x_2)}{G+\exp(x_2)}, x^*\exp(x_3)\right)$ initiated from $(x_1(0), x_2(0), x_3(0))' = x_0$ satisfies:

$$|(x_1(t), x_2(t), x_3(t))| \leq a(|x_0|) \qquad (3.8)$$

Since $z(t) = x(t)$ for all $t \geq r$, it follows from (3.3), (3.5), (3.7) and (3.8) that the closed-loop system (3.1) with (3.2) and $D = k(s, z)$ satisfies the following estimates:

$$\left|\ln\left(\dfrac{x(t)}{x^*}\right)\right| + \left|\ln\left(\dfrac{(s_{in}-s^*)s(t)}{s^*(s_{in}-s(t))}\right)\right| + \left|\ln\left(\dfrac{z(t)}{x^*}\right)\right| \leq 3\sigma\left(a\left(\left|\ln\left(\dfrac{x(0)}{x^*}\right)\right| + \left|\ln\left(\dfrac{z(0)}{x^*}\right)\right| + \left|\ln\left(\dfrac{(s_{in}-s^*)s(0)}{s^*(s_{in}-s(0))}\right)\right|\right), t-r\right), \text{ for all } t \geq r \quad (3.9)$$

$$\left|\ln\left(\dfrac{x(t)}{x^*}\right)\right| + \left|\ln\left(\dfrac{(s_{in}-s^*)s(t)}{s^*(s_{in}-s(t))}\right)\right| + \left|\ln\left(\dfrac{z(t)}{x^*}\right)\right| \leq 3a\left(\left|\ln\left(\dfrac{x(0)}{x^*}\right)\right| + \left|\ln\left(\dfrac{z(0)}{x^*}\right)\right| + \left|\ln\left(\dfrac{(s_{in}-s^*)s(0)}{s^*(s_{in}-s(0))}\right)\right|\right), \text{ for all } t \in [0, r] \quad (3.10)$$

Inequalities (3.8), (3.10) allow us to conclude that the equilibrium point $(s, x, z) = (s^*, x^*, x^*)$ is globally asymptotically stable for the closed-loop system (3.1) with (3.2) and $D = k(s, z)$. The proof is complete. ◁

**Example 3.2:** In [18] it is shown that the feedback law $D = \mu(s)\dfrac{D^* s^*}{(D^*+b)x^*}\dfrac{x}{s} + L\max(0, s^* - s)$, where $L > 0$ is a constant, guarantees that $(s^*, x^*) \in (0, s_{in}) \times (0, +\infty)$ is globally asymptotically stable for the closed-loop system (3.1) with $g(s) = K\mu(s)$, $s_{in}(t) \equiv s_{in}$, where $K > 0$ is a constant.

Here, we will show that system (3.1) with $g(s) = K\mu(s)$, $s_{in}(t) \equiv s_{in}$, $D = \mu(s)\dfrac{D^* s^*}{(D^* + b)x^*}\dfrac{z}{s} + L \max(0, s^* - s)$, $\dot{z} = (\mu(s) - D - b)z$ is forward complete. Therefore, Proposition 3.1 guarantees that for every $r > 0$, the equilibrium point $(s, x, z) = (s^*, x^*, x^*)$ is globally asymptotically stable for the closed-loop system (3.1) with (3.2) and $D = \mu(s)\dfrac{D^* s^*}{(D^* + b)x^*}\dfrac{z}{s} + L \max(0, s^* - s)$, $g(s) = K\mu(s)$, $s_{in}(t) \equiv s_{in}$.

We consider system (3.1) with $g(s) = K\mu(s)$, $s_{in}(t) \equiv s_{in}$, $D = \mu(s)\dfrac{D^* s^*}{(D^* + b)x^*}\dfrac{z}{s} + L \max(0, s^* - s)$, $\dot{z} = (\mu(s) - D - b)z$. Clearly, for every $(z_0, x_0, s_0) \in (0, +\infty) \times (0, +\infty) \times (0, s_{in})$ there exists $t_{\max} > 0$ such that the solution $(z(t), x(t), s(t)) \in (0, +\infty) \times (0, +\infty) \times (0, s_{in})$ of system (3.1) with $g(s) = K\mu(s)$, $s_{in}(t) \equiv s_{in}$, $D = \mu(s)\dfrac{D^* s^*}{(D^* + b)x^*}\dfrac{z}{s} + L \max(0, s^* - s)$, $\dot{z} = (\mu(s) - D - b)z$ initiated from $(z(0), x(0), s(0)) = (z_0, x_0, s_0)$ exists for all $t \in [0, t_{\max})$.

Using the fact that $\mu(s) \leq \mu_{\max}$ for all $s \geq 0$, it follows that:

$$z(t) = x(t)\frac{z_0}{x_0}, \text{ for all } t \in [0, t_{\max}) \tag{3.11}$$

$$x(t) \leq x_0 \exp((\mu_{\max} - b)t), \text{ for all } t \in [0, t_{\max}) \tag{3.12}$$

Simple manipulations show that $\dot{s} = L\max(0, s^* - s)(s_{in} - s) + K\mu(s)x\dfrac{s_{in}}{cs(s_{in} - s^*)}(-\omega s + s^*)$, where $c = \dfrac{x_0}{z_0}$, $\omega = c + (1 - c)\dfrac{s^*}{s_{in}}$. It follows that $\dot{s} < 0$ for all $s > s^* \max(1, \omega^{-1})$ and $\dot{s} > 0$ for all $s < s^* \min(1, \omega^{-1})$. At this point it should be noticed that $s_{in} > s^* \max(1, \omega^{-1})$. Therefore, the following inequalities hold:

$$\min(s_0, s^*, s^*\omega^{-1}) \leq s(t) \leq \max(s_0, s^*, s^*\omega^{-1}), \text{ for all } t \in [0, t_{\max}) \tag{3.13}$$

Using (3.11), (3.12), (3.13), it follows that the following differential inequality holds for all $t \in [0, t_{\max})$

$$\dot{x} \geq \left(-\mu_{\max}\frac{s^*}{(s_{in} - s^*)\min(s_0, s^*, s^*\omega^{-1})} K z_0 \exp((\mu_{\max} - b)t) - Ls^* - b\right)x$$

which directly implies that the following estimate holds:

$$x(t) \geq \exp\left(-\frac{\mu_{\max} s^* K z_0}{(\mu_{\max} - b)(s_{in} - s^*)\min(s_0, s^*, s^*\omega^{-1})}(\exp((\mu_{\max} - b)t) - 1) - Ls^* t - bt\right)x_0, \text{ for all } t \in [0, t_{\max}) \tag{3.14}$$

Inequalities (3.12), (3.13) and (3.14) in conjunction with (3.11) and a standard contradiction argument show that we must have $t_{\max} = +\infty$. Hence, system (3.1) with $g(s) = K\mu(s)$, $s_{in}(t) \equiv s_{in}$, $D = \mu(s)\dfrac{D^* s^*}{(D^* + b)x^*}\dfrac{z}{s} + L \max(0, s^* - s)$, $\dot{z} = (\mu(s) - D - b)z$ is forward complete. ◁

**3.B.** Coexistence implies absence of strong observability for $n = 2$

Consider system (1.1) for $n = 2$ with constant inputs $D(t) \equiv D$ and $s_{in}(t) \equiv s_{in}$. A coexistence equilibrium point for the chemostat model (1.1) is an equilibrium point $(x_1^*, x_2^*, s^*) \in \text{int}(\Re_+^2) \times (0, s_{in})$ of (1.1) satisfying (see [30])

$$\mu_1(s^*) - b_1 = D = \mu_2(s^*) - b_2 \tag{3.15}$$

$$D(s_{in} - s^*) = g_1(s^*) x_1^* + g_2(s^*) x_2^* \tag{3.16}$$

**Proposition 3.3:** *If system (1.1) with $n = 2$ admits a coexistence equilibrium point then for every $r > 0$, system (1.1) with $n = 2$ is not strongly observable in time $r > 0$.*

**Proof:** By virtue of Definition 2.1 and Definition 2.2, it suffices to show that there exist inputs $(D, s_{in}) \in L^\infty([0, r]; U)$, states $(x_0, s_0) \in \Omega \times A$, $(\xi_0, s_0) \in \Omega \times A$ such that

$$s(t) = \bar{s}(t), \text{ for all } t \in [0, r] \tag{3.17}$$

where $s(t) \in A$ denotes the component of the solution $(x(t), s(t)) \in \Omega \times A$ of (1.1), with initial condition $(x_0, s_0) \in \Omega \times A$ corresponding to input $(D, s_{in}) \in L^\infty_{loc}(\Re_+; U)$ and $\bar{s}(t) \in A$ denotes the component of the solution $(\bar{x}(t), \bar{s}(t)) \in \Omega \times A$ of (1.1), with initial condition $(\xi_0, s_0) \in \Omega \times A$ corresponding to input $(D, s_{in}) \in L^\infty_{loc}(\Re_+; U)$.

Consider the constant inputs $D(t) \equiv D$, $s_{in}(t) \equiv s_{in}$ and notice that for every $\xi = (\xi_1, \xi_2) \in \Omega$ with

$$\xi_2 = \frac{g_1(s^*)}{g_2(s^*)}(x_1^* - \xi_1) + x_2^* \tag{3.18}$$

the component $s(t) \in A$ of the solution $(x(t), s(t)) \in \Omega \times A$ of (1.1), with initial condition $(x^*, s^*) \in \Omega \times A$ corresponding to input $(D, s_{in}) \in L^\infty_{loc}(\Re_+; U)$ and the component $\bar{s}(t) \in A$ of the solution $(\bar{x}(t), \bar{s}(t)) \in \Omega \times A$ of (1.1), with initial condition $(\xi, s^*) \in \Omega \times A$ corresponding to input $(D, s_{in}) \in L^\infty_{loc}(\Re_+; U)$ satisfy

$$s(t) \equiv \bar{s}(t) \equiv s^* \tag{3.19}$$

The proof is complete. ◁

**Remark 3.4:** It should be emphasized that if a coexistence equilibrium point for the chemostat model (1.1) with $n = 2$ exists, then no smooth observer of the form

$$\begin{aligned}
\dot{z}_1(t) &= (\mu_1(s(t)) - D(t) - b_1) z_1(t) + (s(t) - \xi(t)) G_1(s(t), \xi(t), z_1(t), z_2(t), D(t), s_{in}(t)) \\
\dot{z}_2(t) &= (\mu_2(s(t)) - D(t) - b_2) z_2(t) + (s(t) - \xi(t)) G_2(s(t), \xi(t), z_1(t), z_2(t), D(t), s_{in}(t)) \\
\dot{\xi}(t) &= D(t)(s_{in}(t) - s(t)) - g_1(s(t)) z_1(t) - g_2(s(t)) z_2(t) + (s(t) - \xi(t)) G_3(s(t), \xi(t), z_1(t), z_2(t), D(t), s_{in}(t)) \\
z(t) &= (z_1(t), z_2(t))' \in \Omega, \xi(t) \in A
\end{aligned} \tag{3.20}$$

for system (1.1) with $n = 2$ exists, where $G_i : A \times A \times \Omega \times U \to \Re$ ($i = 1, 2, 3$) are locally Lipschitz functions. Indeed, using the argument in the proof of Proposition 3.3, we consider system (1.1) with $n = 2$, $D(t) \equiv D^*$, $s_{in}(t) \equiv s_{in}$ initiated at $(x_1, x_2, s^*)$, where $x = (x_1, x_2) \in \Omega$ satisfies $D^*(s_{in} - s^*) = g_1(s^*) x_1 + g_2(s^*) x_2$ and system (3.20) initiated at $(z_1, z_2, s^*)$, where $z = (z_1, z_2) \in \Omega$ satisfies $D^*(s_{in} - s^*) = g_1(s^*) z_1 + g_2(s^*) z_2$ and $z \neq x$. In this case it holds that $\xi(t) = s(t) \equiv s^*$ and $|z(t) - x(t)| \equiv |z(0) - x(0)|$.

It should be noticed that no smooth observer of the form (3.20) can be designed such that for every $(D, s_{in}) \in L_{loc}^\infty(\Re_+; U)$ we have $\lim_{t \to +\infty} |z(t) - x(t)| = 0$. If special inputs $(D, s_{in}) \in L_{loc}^\infty(\Re_+; U)$ are used then we may be able to obtain an observer for system (1.1). The following subsection shows that this is the case.

**3.C.** Batch Culture for $n = 2$.

The model of a batch culture of $n$ species in competition is given by (1.1) with $D \equiv 0$. For $n = 2$ we obtain the input-free model:

$$\dot{x}_1 = (\mu_1(s) - b_1)x_1 \quad , \quad \dot{x}_2 = (\mu_2(s) - b_2)x_2$$
$$\dot{s} = -g_1(s)x_1 - g_2(s)x_2 \tag{3.21}$$

We will assume constant yield coefficients, i.e., $g_i(s) = \mu_i(s)$ for $i = 1,2$ and Michaelis-Menten kinetics for the specific growth rates, i.e.,

$$\mu_1(s) = \frac{a_1 s}{k_1 + s} \quad ; \quad \mu_2(s) = \frac{a_2 s}{k_2 + s} \tag{3.22}$$

where $a_1, a_2, k_1, k_2 > 0$ are positive constants. Therefore, we obtain the model:

$$\dot{x}_1 = \left(\frac{a_1 s}{k_1 + s} - b_1\right)x_1 \quad , \quad \dot{x}_2 = \left(\frac{a_2 s}{k_2 + s} - b_2\right)x_2$$
$$\dot{s} = -\frac{a_1 s}{k_1 + s} x_1 - \frac{a_2 s}{k_2 + s} x_2 \tag{3.23}$$
$$(x_1, x_2, s) \in \mathrm{int}(\Re_+^3)$$

We are in a position to prove the following result.

**Theorem 3.5:** *For every $r > 0$ the following implication holds:*

System (3.23) is not strongly observable in time $r > 0$ $\Rightarrow$ $a_1 = a_2$, $b_1 = b_2$ and $k_1 = k_2$

**Remark 3.6:** Theorem 3.5 guarantees that if the specific growth rates and the mortality rates of the two species are not identical, then for every $r > 0$, system (3.23) is strongly observable in time $r > 0$. Therefore the hybrid observer (2.6) will be a reduced order dead-beat observer. It is clear that the converse implication of the one provided by Theorem 3.2 automatically holds. Notice that if $a_1 = a_2 = a$, $b_1 = b_2 = b$ and $k_1 = k_2 = k$, then system (3.23) can be decomposed into two subsystems:

---the "observable" subsystem

$$\frac{d}{dt}(x_1 + x_2) = \left(\frac{as}{k+s} - b\right)(x_1 + x_2)$$
$$\dot{s} = -\frac{as}{k+s}(x_1 + x_2)$$

and

---the "unobservable" subsystem

$$\frac{d}{dt}(x_1 - x_2) = \left(\frac{as}{k+s} - b\right)(x_1 - x_2)$$

Therefore the characterization provided by Theorem 3.5 is sharp.

**Proof:** Suppose that system (3.23) is not strongly observable in time $r > 0$. Applying Lemma 2.4 with $D \equiv 0$ and using (3.22) and the fact that $g_i(s) = \mu_i(s)$ for $i = 1, 2$, we guarantee the existence of $(x_0, s_0) \in \text{int}(\Re_+^3)$ such that:

$$(k_2 - k_1)\dot{s} = [a_1 - a_2 + b_2 - b_1]s^2 + [a_1 k_2 - a_2 k_1 + (b_2 - b_1)(k_1 + k_2)]s + k_1 k_2 (b_2 - b_1), \text{ for all } t \in [0, r] \quad (3.24)$$

$$\dot{s}(t) = -\left(x_{2,0} + \frac{\mu_1(s_0)}{\mu_2(s_0)} x_{1,0}\right) \mu_2(s(t)) \exp\left(\int_0^t (\mu_2(s(w)) - b_2) dw\right) < 0, \text{ for all } t \in [0, r] \quad (3.25)$$

where $s(t) > 0$ denotes the component of the solution $(x(t), s(t)) \in \text{int}(\Re_+^3)$ of (3.23) with initial condition $(x_0, s_0) \in \text{int}(\Re_+^3)$.

We next distinguish the following cases:

<u>Case 1:</u> $k_1 = k_2 = k$. From (3.24) and (3.25) we conclude that the polynomial $p(s) = [a_1 - a_2 + b_2 - b_1]s^2 + k[a_1 - a_2 + 2(b_2 - b_1)]s + k^2(b_2 - b_1)$ must be the zero polynomial. It follows that $b_1 = b_2$ and $a_1 = a_2$.

<u>Case 2:</u> $k_1 \neq k_2$. We will show that this case cannot happen because it leads to a contradiction.

In this case, we get from (3.24) and (3.25):

$$\dot{s} = f(s) := \frac{a_1 - a_2 + b_2 - b_1}{k_2 - k_1} s^2 + \frac{a_1 k_2 - a_2 k_1 + (b_2 - b_1)(k_1 + k_2)}{k_2 - k_1} s + k_1 k_2 \frac{b_2 - b_1}{k_2 - k_1}, \text{ for all } t \in [0, r] \quad (3.26)$$

and

$$\frac{\dot{s}}{\mu_2(s)} \exp\left(-\int_0^t (\mu_2(s(w)) - b_2) dw\right) = -A, \text{ for all } t \in [0, r] \quad (3.27)$$

where $A := x_{2,0} + \frac{\mu_1(s_0)}{\mu_2(s_0)} x_{1,0}$. By virtue of (3.27) it follows that for all $t \in [0, r]$:

$$\ddot{s} - \frac{\dot{s}^2}{\mu_2(s)} \mu_2'(s) - \dot{s}(\mu_2(s) - b_2) = 0 \quad (3.28)$$

and by virtue of (3.26) it follows that for all $t \in [0, r]$:

$$\ddot{s} = f'(s) f(s) \quad (3.29)$$

Combining (3.26), (3.28) and (3.29), we obtain for all $t \in [0, r]$:

$$\left[f'(s) - \frac{f(s)}{\mu_2(s)} \mu_2'(s) - (\mu_2(s) - b_2)\right] f(s) = 0$$

or equivalently,

$$p(s) f(s) = 0 \quad (3.30)$$

where

$$p(s) = 2(a_1 - a_2 + b_2 - b_1)s^3 + [(2a_1 - 2a_2 + 3b_2 - b_1)k_2 - b_1(k_1 + k_2)]s^2 \\ + k_2 b_2 (k_2 - k_1) s - (b_2 - b_1) k_1 k_2^2 \quad (3.31)$$

The fact that $\dot{s} < 0$ implies that the polynomial $p(s)f(s)$ must be the zero polynomial. Therefore, $b_1 = b_2 = 0$ and $a_1 = a_2 = a > 0$. However, notice that in this case (3.26) gives $\dot{s} = as$, which contradicts (3.25) and the fact $\dot{s} < 0$.

The proof is complete. ◁

**3.D.** Conditions for strong observability in time $r > 0$ for $n = 2$.

We next provide conditions for strong observability of system (1.1) with $n = 2$ for two different cases.

**Theorem 3.8:** *Consider system (1.1) with $n = 2$, $s_{in}(t) \equiv s_{in}$ and $A = (0, s_{in})$. Assume that there exists a constant $c > 0$ such that one of the following holds:*

**(A1)** $\kappa(s) \neq 0$ and $\dfrac{\mu_2(s) - \mu_1(s) + b_1 - b_2}{\kappa(s)} \leq -c$, for all $s \in (0, s_{in})$

or

**(A2)** $\kappa(s) \neq 0$ and $\dfrac{\mu_2(s) - \mu_1(s) + b_1 - b_2}{\kappa(s)} \geq c$, for all $s \in (0, s_{in})$

*where $\kappa$ is defined by (2.10). Then for every $r \geq c^{-1} s_{in}$, (1.1) is strongly observable in time $r > 0$.*

**Theorem 3.9:** *Consider system (1.1) with $n = 2$ and $A = (0, +\infty)$. Assume that there exist constants $a, c > 0$ such that one of the following holds:*

**(A3)** $\kappa(s) \neq 0$ and $\dfrac{\mu_2(s) - \mu_1(s) + b_1 - b_2}{\kappa(s)} \leq -as^2 - c$, for all $s \in (0, +\infty)$

or

**(A4)** $\kappa(s) \neq 0$ and $\dfrac{\mu_2(s) - \mu_1(s) + b_1 - b_2}{\kappa(s)} \geq as^2 + c$, for all $s \in (0, +\infty)$

*where $\kappa$ is defined by (2.10). Then for every $r \geq c^{-1} + a^{-1}$, (1.1) is strongly observable in time $r > 0$.*

**Remark 3.10:** For the usual case where $g_i(s) = \mu_i(s)$ for $i = 1, 2$ and the specific growth rates satisfy the Michaelis-Menten kinetics (3.22), conditions (A1), (A2), (A3), (A4) are expressed by:

**(A1')** $k_1 \neq k_2$ and $\dfrac{a_1 - a_2 + b_2 - b_1}{k_2 - k_1} s^2 + \dfrac{a_1 k_2 - a_2 k_1 + (b_2 - b_1)(k_1 + k_2)}{k_2 - k_1} s + k_1 k_2 \dfrac{b_2 - b_1}{k_2 - k_1} \leq -c$, for all $s \in (0, s_{in})$

or

**(A2')** $k_1 \neq k_2$ and $\dfrac{a_1 - a_2 + b_2 - b_1}{k_2 - k_1} s^2 + \dfrac{a_1 k_2 - a_2 k_1 + (b_2 - b_1)(k_1 + k_2)}{k_2 - k_1} s + k_1 k_2 \dfrac{b_2 - b_1}{k_2 - k_1} \geq c$, for all $s \in (0, s_{in})$

and

**(A3')** $k_1 \neq k_2$ and $\dfrac{a_1 - a_2 + b_2 - b_1}{k_2 - k_1} s^2 + \dfrac{a_1 k_2 - a_2 k_1 + (b_2 - b_1)(k_1 + k_2)}{k_2 - k_1} s + k_1 k_2 \dfrac{b_2 - b_1}{k_2 - k_1} \leq -as^2 - c$, for all $s \in (0, +\infty)$

or

**(A4')** $k_1 \neq k_2$ and $\dfrac{a_1 - a_2 + b_2 - b_1}{k_2 - k_1} s^2 + \dfrac{a_1 k_2 - a_2 k_1 + (b_2 - b_1)(k_1 + k_2)}{k_2 - k_1} s + k_1 k_2 \dfrac{b_2 - b_1}{k_2 - k_1} \geq as^2 + c$, for all $s \in (0, +\infty)$

For Theorem 3.8, if in addition we have $D(t) \in [0, D_{max}]$, where $D_{max} \geq 0$, then condition (A2) can obtain the following, less demanding form:

**(A2'')** $\kappa(s) \neq 0$ for all $s \in (0, s_{in})$ and $\dfrac{\mu_2(s) - \mu_1(s) + b_1 - b_2}{\kappa(s)} \geq c$ for all $s \in (0, s_{in})$ with $\dfrac{\mu_2(s) + b_1 - \mu_1(s) - b_2}{\kappa(s)(s_{in} - s)} < D_{max}$.

**Proof of Theorem 3.8 and Theorem 3.9:** Employing Lemma 2.4, we can guarantee the existence of $(D, s_{in}) \in L^\infty([0, r]; U)$ and $(x_0, s_0) \in \text{int}(\Re_+^2) \times A$ such that:

$$\dot{s}(t) = \dfrac{\mu_2(s(t)) - \mu_1(s(t)) + b_1 - b_2}{\kappa(s(t))}, \text{ for almost all } t \in [0, r] \quad (3.32)$$

Notice that (3.32) can only be satisfied $\dot{s}(t) = \dfrac{\mu_2(s(t)) - \mu_1(s(t)) + b_1 - b_2}{\kappa(s(t))} < D(t)(s_{in}(t) - s(t))$ and consequently, for the case of Theorem 3.8 with $D(t) \in [0, D_{max}]$, we can guarantee that $\dfrac{\mu_2(s(t)) + b_1 - \mu_1(s(t)) - b_2}{\kappa(s(t))(s_{in} - s(t))} < D_{max}$ for all $t \in [0, r]$. It is direct to verify that if $s_0 \in A = (0, s_{in})$, $r \geq c^{-1} s_{in}$ and hypothesis (A1) or hypothesis (A2) holds then the solution of (3.32) cannot satisfy $s(t) \in A = (0, s_{in})$ for all $t \in [0, r]$, a contradiction.

For Theorem 3.9, if hypothesis (A4) holds then by using the comparison lemma in [21], we can guarantee that the solution of (3.32) presents a finite escape time in the interval $[0, a^{-1} + c^{-1}]$. On the other hand, if hypothesis (A3) holds then by employing the comparison lemma in [21], we can guarantee that the solution of (3.32) cannot satisfy $s(t)$ for all $t \in [0, a^{-1} + c^{-1}]$.

Details are left to the reader. ◁

## 4. Concluding Remarks

In this work, we have applied recent results for the observability of systems which are linear with respect to the unmeasured state components in [20] in order to design hybrid dead-beat reduced-order observers for the chemostat with microbial species in competition, i.e., system (1.1). We have assumed that the measured output is the concentration of the nutrient and we are interested in the estimation of the size of the populations of the competing microbial species. We have showed that the chemostat with one species is strongly observable in time $r > 0$ for arbitrary $r > 0$. The design of a reduced-order hybrid dead-beat observer for the case with one species allowed us to show that the dynamic output feedback stabilization problem for (1.1) with $n = 1$ can be solved with the combination of a static state feedback stabilizer and the proposed dead-beat hybrid reduced order observer.

Furthermore, new relationships between coexistence and strong observability are established for a chemostat with two microbial species (Proposition 3.3). The proposed dead-beat reduced-order observer can be used for system (1.1) with $n = 2$, $D \equiv 0$ (batch culture) and Michaelis-Menten kinetics for the specific growth rates: the batch culture is strongly observable in time $r > 0$ for arbitrary $r > 0$ if and only if the specific growth rates and the mortality rates of the two microbial species are not identical. The result (Theorem 3.5) is important because batch cultures of microbial species are used only for finite time and the proposed dead-beat hybrid reduced order observer can provide exact estimates in very short times. Finally, a set of sufficient conditions (which do not allow coexistence) that can guarantee strong observability of system (1.1) with $n = 2$ in time $r > 0$ is provided.

The obtained results can be extended in the same spirit to the case of system (1.1) with $n \geq 3$. Moreover, the obtained results can be used for the study of the observer design problem of system (1.1) under periodic inputs $(D, s_{in}) \in L^{\infty}_{loc}(\Re_+; U)$. This will be the subject of future research.

**Acknowledgements:** This work has been supported in part by the NSF grants DMS-0504462 and DMS-0906659.

# References

[1] Angeli, D. and E.D. Sontag, "Forward Completeness, Unbounded Observability and Their Lyapunov Characterizations", *Systems and Control Letters*, 38(4-5), 1999, 209-217.

[2] Angeli, D., P. De Leenheer and E. D. Sontag, "A Small-Gain Theorem for Almost Global Convergence of Monotone Systems", *Systems and Control Letters*, 52(5), 2004, 407-414.

[3] Bastin, G. and D. Dochain, *On-line Estimation and Adaptive Control of Bioreactors*, Elsevier, Amsterdam, 1990.

[4] Bernard, O., G. Sallet and A. Sciandra, "Nonlinear Observers for a Class of Biological Systems Application to Validation of Phytoplanktonic Growth Model", *IEEE Transactions on Automatic Control*, 43, 1998, 1056-1065.

[5] Bogaerts, P., "A Hybrid Asymptotic-Kalman Observer for Bioprocesses", *Bioprocess Engineering*, 20(3), 1999, 101-113.

[6] De Leenheer, P. and H. L. Smith, "Feedback Control for Chemostat Models", *Journal of Mathematical Biology*, 46, 2003, 48-70.

[7] De Leenheer, P., D. Angeli and E. D. Sontag, "Crowding Effects Promote Coexistence in the Chemostat", *Journal of Mathematical Analysis and Applications*, 319, 2006, 48-60.

[8] Dumont, M., A. Rapaport, J. Harmand, B. Benyahia and J.-J. Godon, "Observers for Microbial Ecology- How Including Molecular Data Into Bioprocess Modeling?", *Proceedings of the 16th Mediterranean Conference on Control and Automation*, France, 2008, 1381-1386.

[9] Farza, M., K. Busawon and H. Hammouri, "Simple Nonlinear Observers for On-Line Estimation of Kinetic Rates in Bioreactors", *Automatica*, 34, 1998, 301-318.

[10] Gauthier, J. P., H. Hammouri and S. Othman, "A Simple Observer for Nonlinear Systems. Application to Bioreactors", *IEEE Transactions on Automatic Control*, 37, 1992, 875-880.

[11] Gouze, J. L., A. Rapaport and M. Z. Hadj-Sadok, "Interval Observers for Uncertain Biological Systems", *Ecological Modelling*, 133, 2000, 45-56.

[12] Gouze, J. L. and G. Robledo, "Feedback Control for Nonmonotone Competition Models in the Chemostat", *Nonlinear Analysis: Real World Applications*, 6, 2005, 671-690.

[13] Gouze, J. L. and G. Robledo, "Robust Control for an Uncertain Chemostat Model", *International Journal of Robust and Nonlinear Control*, 16(3), 2006, 133-155.

[14] Harmard, J., A. Rapaport and F. Mazenc, "Output Tracking of Continuous Bioreactors Through Recirculation and By-Pass", *Automatica*, 42, 2006, 1025-1032.

[15] Karafyllis, I., "Non-Uniform in Time Robust Global Asymptotic Output Stability", *Systems and Control Letters*, 54(3), 2005, 181-193.

[16] Karafyllis, I., C. Kravaris, L. Syrou and G. Lyberatos, "A Vector Lyapunov Function Characterization of Input-to-State Stability With Application to Robust Global Stabilization of the Chemostat", *European Journal of Control*, 14(1), 2008, 47-61.

[17] Karafyllis, I. and C. Kravaris, "Global Stability Results for Systems Under Sampled-Data Control", *International Journal of Robust and Nonlinear Control*, 19(10), 2009, 1105-1128.

[18] Karafyllis, I., C. Kravaris and N. Kalogerakis, "Relaxed Lyapunov Criteria for Robust Global Stabilization of Nonlinear Systems", *International Journal of Control*, 82(11), 2009, 2077-2094.

[19] Karafyllis, I. and Z.-P. Jiang, "New Results in Trajectory-Based Small-Gain With Application to the Stabilization of a Chemostat", submitted to the *International Journal of Robust and Nonlinear Control*. (see also http://arxiv.org/abs/1002.4489 )

[20] Karafyllis, I. and Z.-P. Jiang, "Hybrid Dead-Beat Observers for a Class of Nonlinear Systems", submitted to *Systems and Control Letters*. (see also http://arxiv.org/abs/1005.3984 )

[21] Khalil, H. K., *Nonlinear Systems*, 2nd Edition, Prentice-Hall, 1996.

[22] Lemesle, V. and J. L. Gouze, "Hybrid Bounded Error Observers for Uncertain Bioreactor Models", *Bioprocess Biosystems Engineering*, 27, 2005, 311-318.

[23] Mailleret, L., J. L. Gouze and O. Bernard, "Global Stabilization of a Class of Partially Known Nonnegative Systems", *Automatica*, 44(8), 2008, 2128-2134.

[24] Malisoff, M. and F. Mazenc, *Constructions of Strict Lyapunov Functions*, Springer Verlag, London, 2009.

[25] Mazenc, F., M. Malisoff and P. De Leenheer, "On the Stability of Periodic Solutions in the Perturbed Chemostat", *Mathematical Biosciences and Engineering*, 4(2), 2007, 319-338.

[26] Mazenc, F., M. Malisoff and J. Harmand, "Further Results on Stabilization of Periodic Trajectories for a Chemostat With Two Species", *IEEE Transactions on Automatic Control*, 53(1), 2008, 66-74.
[27] Rapaport, A. and J. Harmand, "Robust Regulation of a Class of Partially Observed Nonlinear Continuous Bioreactors", *Journal of Process Control,* 12, 2002, 291-302.
[28] Rapaport, A., J. Harmand and F. Mazenc, "Coexistence in the Design of a Series of Two Chemostats", *Nonlinear Analysis: Real World Applications*, 9, 2008, 1052-1067.
[29] Smith, H., L., *Monotone Dynamical Systems An Introduction to the Theory of Competitive and Cooperative Systems*, Mathematical Surveys and Monographs, Volume 41, AMS, Providence, Rhode Island, 1994.
[30] Smith, H. and P. Waltman, *The Theory of the Chemostat. Dynamics of Microbial Competition*, Cambridge Studies in Mathematical Biology, 13, Cambridge University Press: Cambridge, 1995.
[31] Sontag, E.D., *Mathematical Control Theory*, 2$^{nd}$ Edition, Springer-Verlag, New York, 1998.
[32] Zhu, L. and X. Huang, "Multiple Limit Cycles in a Continuous Culture Vessel With Variable Yield", *Nonlinear Analysis*, 64, 2006, 887-894.
[33] Zhu, L., X. Huang and H. Su, "Bifurcation for a Functional Yield Chemostat When One Competitor Produces a Toxin", *Journal of Mathematical Analysis and Applications*, 329, 2007, 891-903.